\documentclass[12pt,reqno]{amsart}
\usepackage{amsmath,amssymb,amscd,latexsym,amsthm,mathrsfs}
\usepackage[unicode]{hyperref}
\usepackage{hypbmsec}
\ifx\pdfoutput\undefined\else
\hypersetup{
colorlinks=true,
linkcolor=blue,
citecolor=blue,
urlcolor=blue,
filecolor=blue,
bookmarksnumbered=true,
pdfstartview=FitH,
pdfhighlight=/N
}
\fi

\newtheorem{theorem}{Theorem}
\newtheorem{claim}{Claim}
\newtheorem{lemma}{Lemma}

\theoremstyle{remark}

\newtheorem{definition}{Definition}

\newtheorem{observation}{Observation}

\def\sech{{\operatorname{sech}}}

\begin{document}

\hypersetup{pdfauthor={Mathew Rogers},%
pdftitle={Odd zeta values}}

\title{Variations of the Ramanujan polynomials and remarks on $\zeta(2j+1)/\pi^{2j+1}$}

\author{Matilde Lal\'{i}n}
\address{D\'epartement de math\'ematiques et de statistique , Universit\'e de Montr\'eal Montreal, QC, Canada}
\email{mlalin@dms.umontreal.ca}
\thanks{M. L. is supported by NSERC Discovery Grant 355412-2008 and a start-up grant from Universit\'e de Montr\'eal. M. R. is supported by NSF award DMS-0803107}

\author{Mathew Rogers}
\address{Department of Mathematics, University of Illinois, Urbana, IL, USA}
\email{mdrogers@illinois.edu}

\date{June 6, 2011}
\subjclass[2000]{Primary ; Secondary }

\begin{abstract}
We observe that five polynomial families have all of their
zeros on the unit circle.  We prove the statements explicitly for
four of the polynomial families.  The polynomials have
coefficients which involve Bernoulli numbers, Euler numbers, and
the odd values of the Riemann zeta function. These polynomials are
closely related to the Ramanujan polynomials, which were recently
introduced by Murty, Smyth and Wang \cite{MSW}. Our proofs rely
upon theorems of Schinzel \cite{Sc}, and Lakatos and Losonczi
\cite{LL} and some generalizations.
\end{abstract}

\maketitle

\section{Introduction}

In a recent paper, Murty, Smyth and Wang considered the
\textit{Ramanujan polynomials} \cite{MSW}.  They were defined by Gun, Murty and Rath \cite{GMR} using
\begin{equation}\label{Ramanujan polynomials}
R_{2k+1}(z):=\sum_{j=0}^{k-1}\frac{B_{2j}B_{2k+2-2j}}{(2j)!(2k+2-2j)!}z^{2j},
\end{equation}
where $B_j$ is the $j$th Bernoulli number.  Among other fascinating
results, Murty, Smyth and Wang showed that $R_{2k+1}(z)$ has all of its non-real
zeros on the unit circle.  The purpose of this paper is to study
some variants of $R_{2k+1}(z)$, which also have many zeros on the
unit circle.
\begin{claim}\label{main theorem}  Let $B_j$ denote the Bernoulli numbers, and let
$E_j$ denote the Euler numbers.  Suppose that $k\ge 2$. The
following polynomials have all of their non-zero roots on the
unit circle:
\begin{align}
P_k(z):=&\frac{(2\pi)^{2k-1}}{(2k)!}\sum_{j=0}^{k}(-1)^j
B_{2j}B_{2k-2j}{2k\choose 2j}z^{2j}\\
&+\zeta(2k-1)\left(z^{2k-1}+(-1)^k z\right),\notag\\
Q_k(z):=&\left(2^{2k}+1\right)P_{k}(z)-2^{2k}P_{k}\left(z/2\right)-P_{k}(2z),\\
Y_{k}(z):=&\frac{\pi}{2^{2k}}\left(Q_{k}(i \sqrt{z})+Q_{k}(-i \sqrt{z})\right)\\
W_k(z):=&\left(2^{2k-1}+2\right)P_{k}(z)-2^{2k}P_{k}\left(z/2\right)-P_{k}(2z),\\
S_k(z):=&\sum_{j=0}^{k}E_{2j}E_{2k-2j}{2k\choose 2j}z^{j}.
\end{align}
\end{claim}
We will offer a general proof of Claim \ref{main theorem} for
$Q_{k}(z)$, $Y_{k}(z)$, $W_k(z)$, and $S_{k}(z)$.  It seems that
$P_{k}(z)$ is more difficult to handle.  In Section \ref{Section:
partial Pk results} we offer several partial results concerning
$P_k(z)$.

An important secondary goal of this work, is to highlight a
connection with the odd values of the Riemann zeta function.
Recall that the Riemann zeta function is defined by
\begin{equation*}
\zeta(s):=\sum_{n=1}^{\infty}\frac{1}{n^s}.
\end{equation*}
While it is a classical fact that $\zeta(2j)/\pi^{2j}$ is rational
when $j\ge 1$, very little is known about the arithmetic nature of
$\zeta(2j+1)$.  The only theorems in this direction are celebrated
irrationality results. For instance, Ap\'{e}ry showed that
$\zeta(3)$ is irrational \cite{Ap, vP}, Rivoal proved that infinitely
many odd zeta values are irrational \cite{Ri}, and Zudilin
established that at least one element of the set $\{\zeta(5),
\zeta(7), \zeta(9), \zeta(11)\}$ is also irrational \cite{Zd}.

Therefore, let us briefly consider the identity which gave birth to
the Ramanujan polynomials \cite[~pg. 276]{Be2}.  The formula can be
written as
\begin{equation}\label{zeta(3) series}
\begin{split}
\frac{(2\pi)^{2k-1}}{2(2k)!}&\sum_{j=0}^{k}(-1)^j
B_{2j}B_{2k-2j}{2k\choose
2j}z^{2j}+\frac{\zeta(2k-1)}{2}\left((-1)^k
z+z^{2k-1}\right)\\
&=-\sum_{n=1}^{\infty}\frac{1}{n^{2k-1}}\frac{z^{2k-1}}{e^{2\pi n
/z}-1}+(-1)^{k+1}\sum_{n=1}^{\infty}\frac{1}{n^{2k-1}}\frac{z}{e^{2\pi
n z}-1}.
\end{split}
\end{equation}
This identity holds whenever $z\not\in i\mathbb{Q}$. The
restriction is necessary to ensure that both infinite series
converge. Not surprinsigly, this formula is also mentioned in works dealing with the irrationality of $\zeta(3)$ \cite{vP}. Notice that Gun, Murty and Rath used \eqref{zeta(3)
series} to express odd zeta values in terms of Eichler integrals
\cite{GMR}. Now consider the case when $k=2$.  A brief numerical
calculation shows that the polynomial obtained from the left-hand
side
\begin{equation}\label{first zeta 3 polynomial}
z^4+5z^2+1-\frac{90\zeta(3)}{\pi^3}\left(z^3+z\right)=0,
\end{equation}
has all of its zeros on the unit circle.  
Notice that if we truncate the right-hand side of \eqref{zeta(3)
series}, then we can nearly obtain a closed form for $\zeta(3)$:
\begin{equation}\label{zeta 3 approx 1}
\zeta(3)\approx\left(\frac{z^4+5z^2+1}{z^3+z}\right)\frac{\pi^3}{90},
\end{equation}
where
\begin{equation*}
0=\frac{z}{e^{2\pi/z}-1}+\frac{z^{-1}}{e^{2\pi z}-1}.
\end{equation*}
This approximation gives six decimal places of numerical accuracy.
The accuracy can be increased by including higher order terms in
the truncation. It would be extremely interesting if this idea
could be used to say something about the irrationality of
$\zeta(3)/\pi^{3}$. Unfortunately such a theorem is well beyond
the scope of this paper.  We will settle for the more modest goal
of studying the polynomial families listed in Claim \ref{main
theorem}.


\section{$S_{k}(z)$, $Y_{k}(z)$ and the theorems of Schinzel, Lakatos and Losonczi}\label{section: First analytic proof}


In this section we prove that $S_{k}(z)$ and $Y_{k}(z)$ have all
of their non-zero roots on the unit circle. Our proofs follow
from applying the theorems of Schinzel \cite{Sc}, Lakatos and
Losonczi \cite{LL}, and Lakatos \cite{L}. Lakatos proved that any
reciprocal polynomial $\sum_{j=0}^{k}A_j z^j$, with real-valued
coefficients, which satisfies
\begin{equation}\label{Lakatos condition 1}
|A_k|\ge\sum_{j=0}^{k}\left| A_{j}-A_{k}\right|,
\end{equation}
must have all of its zeros on the unit circle.  If the inequality
is strict then the polynomial has only simple zeros.  Equation
\eqref{Lakatos condition 1} is a very strong restriction. There
have been a number of recent improvements to \eqref{Lakatos
condition 1} with a similar flavor (see \cite{Sc} and \cite{LL}).
Schinzel proved that any self-inversive polynomial which satisfies
\begin{equation}\label{Lakatos condition 2}
|A_k|\ge\inf_{\substack{c,d\in\mathbb{C}\\|d|=1}}\sum_{j=0}^{k}\left|c
A_{j}-d^{k-j} A_{k}\right|,
\end{equation}
must have all of its zeros on the unit circle \cite{Sc}.
Self-inversive polynomials have complex-valued coefficients which
satisfy $A_j=\epsilon\overline{A_{k-j}}$, for some fixed
$|\epsilon|=1$. Notice that the class of self-inversive
polynomials includes both reciprocal and anti-reciprocal
polynomials.  In Theorems \ref{S theorem} and \ref{Theorem on qk}
we apply Schinzel's theorem with $d=1$.

\begin{theorem}\label{S theorem} Suppose that $k\ge 1$.  Then all of the zeros of \[S_k(z)=\sum_{j=0}^k E_{2j} E_{2k-2j} \binom{2k}{2j}
z^{j}\] lie on the unit circle.  Furthermore, all of the zeros are
simple.
\end{theorem}
\begin{proof}
With \eqref{Lakatos condition 2} in mind, let us begin by setting
\[A_j:= E_{2j} E_{2k-2j} \binom{2k}{2j}.\]
The sign of
$E_{2n}$ alternates with respect to $n$.  This implies that all of
the coefficients of $S_k(z)$ have sign $(-1)^k$.

Our proof consists of three main steps. First we remove the
absolute values from the sum \eqref{Lakatos condition 2}. This is
easily accomplished by showing that $(-1)^k (c A_j-A_{k-2})>0$ for
$c=\frac{\pi}{4(1+3^{-1-2k})}$. Next we evaluate $\sum_{j=0}^{k}
A_j$ explicitly, and finally we deduce the desired upper bound.

In order to remove the absolute value
signs from \eqref{Lakatos condition 2}, we need to
demonstrate that $(-1)^k (c A_j-A_{k-2})>0$.  Using
the following bound for Euler numbers \cite[~pg. 805]{AbSt}:
\[\frac{4^{n+1}(2n)!}{\pi^{2n+1}}>|E_{2n}|>\frac{4^{n+1}(2n)!}{\pi^{2n+1}(1+3^{-1-2n})},\]
and the fact that $E_{0}=1$, leads to
\begin{align*}
\frac{4(1+3^{-1-2k})}{\pi}|A_k|&=
\frac{4(1+3^{-1-2k})}{\pi}|E_{2k}E_{0}|\\
&> \frac{4^{k+2}(2k)!}{\pi^{2k+2}}\\
&> |E_{2j} E_{2k-2j}| \binom{2k}{2j}=|A_j|.
\end{align*}
The absolute values can be removed because both $A_k$ and $A_j$ have
sign $(-1)^k$:
\[\frac{(-1)^k4(1+3^{-1-2k})}{\pi}A_k > (-1)^k A_j>0.\]
If we take \[c=\frac{\pi}{4(1+3^{-1-2k})},\] then the previous
inequality implies $(-1)^kA_k-c(-1)^kA_j>0.$

Let $E_n(z)$ denote the classical Euler polynomials, and recall a
standard convolution identity \cite{Di}:
\[\sum_{j=0}^n\binom{n}{j}E_j(v)E_{n-j}(w)=2(1-w-v)E_n(v+w)+2E_{n+1}(v+w).\]
Set $v=w=\frac{1}{2}$, and then use
$E_n=2^nE_n\left(\frac{1}{2}\right)$, to obtain an expression for
$S_{k}(1)$.  We have
\[|S_k(1)|=\left|2^{2k+1}E_{2k+1}(1)\right|=\frac{2^{2k+1}(2^{2k+2}-1)| B_{2k+2}|}{k+1}.\]
The evaluation of $E_{2k+1}(1)$ follows from \cite[~pg. 805]{AbSt}.
  Thus
\begin{align*}
\sum_{j=0}^k |cA_j-A_k|=& (-1)^k \sum_{j=0}^k (A_k-cA_j)\\
=&(k+1)(-1)^k A_k- c (-1)^k S_{k}(1)\\
=&(k+1)|A_k|- c |S_{k}(1)|\\
=& (k+1)|A_k|-\frac{\pi
2^{2k-1}(2^{2k+2}-1)|B_{2k+2}|}{(k+1)(1+3^{-1-2k})}.
\end{align*}
To finish the verification of \eqref{Lakatos condition 2}, we need
to show that the last expression is bounded from above by $|A_k|$.
This is equivalent to showing that
\begin{equation}\label{Sk proof final inequality}
\frac{\pi 2^{2k-1}(2^{2k+2}-1)|B_{2k+2}|}{(k+1)(1+3^{-1-2k})}\geq
k|A_{k}|=k|E_{2k}|.
\end{equation}
We will resort to an inequality for Bernoulli numbers
\cite[~pg. 805]{AbSt}:
\begin{equation}\label{ineq:bernoulli}
\frac{2(2n)!}{(2\pi)^{2n}}<|B_{2n}|<\frac{2(2n)!}{(2\pi)^{2n}(1-2^{1-2n})}.
\end{equation}
Thus we find
\begin{align}
\frac{\pi 2^{2k-1}(2^{2k+2}-1)|B_{2k+2}|}{(k+1)(1+3^{-1-2k})}>&
\frac{\pi
2^{2k-1}(2^{2k+2}-1)}{(k+1)(1+3^{-1-2k})}\frac{2(2k+2)!}{(2\pi)^{2k+2}}\notag\\
=&\frac{2^{2k+1}(1-2^{-2-2k})}{(1+3^{-1-2k})}\frac{(2k+1)!}{\pi^{2k+1}}\label{Berninequality}.
\end{align}
On the other hand, we have already used the fact that Euler numbers
are bounded from above by
\begin{equation}\label{Eulerinequality}
k|E_{2k}| < k \frac{4^{k+1}(2k)!}{\pi^{2k+1}}.
\end{equation}
Substituting \eqref{Berninequality} and \eqref{Eulerinequality} into
\eqref{Sk proof final inequality}, reduces the inequality to
\[\frac{2k+1}{2k}>\frac{1+3^{-1-2k}}{1-2^{-2-2k}}.\]
This final inequality is easily verified with elementary calculus
for $k\ge 1$. Since the inequality is strict, we conclude
immediately that $S_{k}(z)$ has only simple zeros which all lie on
the unit circle.
\end{proof}

We have proved that $S_{k}(z)$ has all of its zeros on the unit
circle.  Perhaps it is interesting to note that $S_{k}(z)$ satisfies
\begin{equation*}
\begin{split}
 \frac{(\pi/2)^{2k+1}}{2(2k)!}S_{k}(-
z^2)=&z^{2k}\sum_{n=1}^{\infty}\frac{\chi_{-4}(n)\sech\left(\pi
n/2z\right)}{n^{2k+1}}\\
&+(-1)^k\sum_{n=1}^{\infty}\frac{\chi_{-4}(n)\sech\left(\pi n
z/2\right)}{n^{2k+1}},
\end{split}
\end{equation*}
where $\chi_{-4}(n)$ is the non-principle character mod $4$. This formula appears immediately
after equation \eqref{zeta(3) series} in Ramanujan's notebook \cite[~pg. 276]{Be2}. As a
result it is easy to approximate the zeros of $S_{k}(z)$ by the
zeros of exponential polynomials.  It remains to be seen whether or
not there are any interesting applications for this observation.

To illustrate this method a second time, we prove that the
polynomial $Y_k(z)$ has all of its non-zero roots on the unit
circle.  Notice that $Y_k(z)$ is a close analogue of $S_k(z)$,
except that it involves Bernoulli numbers rather than Euler
numbers.

\begin{theorem}\label{Theorem on qk}  Suppose that $k\ge 2$. The polynomial
\begin{align*}
Y_{k}(z)=&\frac{\pi}{2^{2k}}\left(Q_{k}(i \sqrt{z})+Q_{k}(-i z)\right)\\
=&\frac{\pi^{2k}}{(2k)!}\sum_{j=0}^k B_{2j}B_{2k-2j}
(2^{2j}-1)(2^{2k-2j}-1)\binom{2k}{2j} z^{j},
\end{align*}
has all of its non-zero roots on the unit circle.  Furthermore,
all of the zeros are simple.
\end{theorem}
\begin{proof}
Observe that $Y_{k}(z)$ has degree $k-1$, since the coefficients
of $z^{k}$ and $z^{0}$ are identically zero. We prove that
$Y_{k}(z)/z$ satisfies the hypothesis of Schinzel's theorem
\cite{Sc}.  If we eliminate the trivial factor of $z$, then we
obtain a polynomial of the form
\begin{equation*}
\frac{Y_{k}(z)}{z}=\sum_{j=0}^{k-2}A_{j} z^{j},
\end{equation*}
where
\begin{equation*}\label{a(q) inversion on a ray}
A_j=
\frac{(2\pi)^{2k}}{(2k)!}\binom{2k}{2j+2}(1-2^{-2-2j})(1-2^{2-2k+2j})B_{2j+2}
B_{2k-2j-2}.
\end{equation*}
Notice that $Y_{k}(z)/z$ is reciprocal, since $A_{k-2-j}=A_{j}$.
By elementary properties of Bernoulli numbers, the sign of $A_j$
is $(-1)^k$ for all $j$.

Schinzel's theorem can be applied if the following inequality holds:
\begin{align}\label{qk proof Schinzel inequality}
|A_{k-2}|\geq& \sum_{j=0}^{k-2}|c A_j -A_{k-2}|,
\end{align}
for some $c\in\mathbb{C}$.  We prove that \eqref{qk proof Schinzel
inequality} holds when
$c=\frac{\pi^2(1-2^{2-2k})}{8(1-2^{3-2k})}$. Our proof follows the
same three steps as in the case of $S_k(z)$.

In order
to remove the absolute value signs from \eqref{qk proof Schinzel
inequality}, we need to demonstrate that $(-1)^k (c
A_j-A_{k-2})>0$.  We demonstrate this by comparing an upper bound
on $(-1)^k A_{k-2}$, with a lower bound on $(-1)^k A_j$. The lower
bound on $|A_j|$ is a consequence of an inequality from
\cite{DAn}:
\[|B_{2n}|>\frac{2(2n)!}{(2\pi)^{2n}(1-2^{-2n})}.\]
In particular we find
\begin{equation}\label{Thm2 Aj lower bound}
(-1)^k A_j=|A_j|>4.
\end{equation}
By \eqref{ineq:bernoulli}, we find an upper bound for $|A_{k-2}|$:
\begin{equation}\label{Thm2 Ak upper bound}
(-1)^k
A_{k-2}=|A_{k-2}|<\frac{\pi^2}{2}\frac{1-2^{2-2k}}{1-2^{3-2k}}.
\end{equation}
 Comparing \eqref{Thm2 Ak upper bound} and
\eqref{Thm2 Aj lower bound}, allows us to easily conclude
\begin{equation}\label{Thm2 Aj Ak inequality}
(-1)^k( c A_j-A_{k-2})>0,
\end{equation}
whenever $k>2$.

Since we have proved \eqref{Thm2 Aj Ak inequality}, Schinzel's sum
immediately reduces to
\begin{equation}\label{qk thm midway} \sum_{j=0}^{k-2}|c A_j
-A_{k-2}|=-(k-1)(-1)^k A_{k-2}+(-1)^k c\sum_{j=0}^{k-2}A_j.
\end{equation}
Now we simplify the remaining sum.  Let $B_{j}(z)$ denote the
usual Bernoulli polynomials.  By standard evaluations of Bernoulli
polynomials \cite[~pg. 805]{AbSt}, we have
\begin{equation*}
A_j=\frac{(2\pi)^{2k}}{4(2k)!}\binom{2k}{2j+2}
\left(B_{2j+2}\left(\frac12\right)-B_{2j+2}(0)\right)\left(
B_{2k-2j-2}\left(\frac12\right)-B_{2k-2j-2}(0)\right).
\end{equation*}
Next we use a well known convolution identity for Bernoulli
polynomials \cite{Di}:
\begin{equation}\label{bernoulliconvolution}
\sum_{j=0}^n\binom{n}{j}B_j(v)B_{n-j}(w)=n(w+v-1)B_{n-1}(v+w)-(n-1)B_{n}(v+w).
\end{equation}
Considering all of the cases where  $(v,w)\in\{\left(\frac12,\frac12
\right),\left(0,\frac12 \right),\left(\frac12,0
\right),\left(0,0\right)\}$, leads to
\begin{align}
\sum_{j=0}^{k-2}A_j&=\frac{(2\pi)^{2k}}{(2k)!}\frac{(2k-1)}{2}\left(B_{2k}\left(\frac12\right)-B_{2k}(0)\right)\notag\\
&= -\frac{(2\pi)^{2k}}{(2k)!}(2k-1)(1-2^{-2k})B_{2k}.\label{qkAjsum}
\end{align}
Substituting \eqref{qkAjsum} into \eqref{qk proof Schinzel
inequality}, leads to a closed form expression for the sum we are
interested in:
\begin{equation*}
\begin{split}
\sum_{j=0}^{k-2}|c A_j -A_{k-2}|=&-(k-1)(-1)^k A_{k-2}-(-1)^k
c\frac{(2\pi)^{2k}}{(2k)!}(2k-1)(1-2^{-2k})B_{2k}.
\end{split}
\end{equation*}
The proof can be completed by showing that this last expression is
bounded from above by $|A_{k-2}|$ or
\[\frac{\pi^2(1-2^{2-2k})}{8(1-2^{3-2k})}\frac{(2\pi)^{2k}}{(2k)!}(2k-1)(1-2^{-2k})|B_{2k}|< k |A_{k-2}|.\]
Employing
\begin{equation}\label{zeta}
\zeta(2n)=\frac{(-1)^{n+1}B_{2n}(2\pi)^{2n}}{2(2n)!}
\end{equation}
reduces the desired inequality to
\begin{equation*}
\frac{1}{2}\left(2k-1\right)(1-2^{-2k})\zeta(2k)<k
(1-2^{3-2k})\zeta(2k-2).
\end{equation*}
It is elementary to show that this inequality holds for $k>1$.
\end{proof}

\section{Generalizing the criteria to other
families}\label{Section on Wk and Qk}

Conditions such as \eqref{Lakatos condition 2} appear to be too
restrictive to apply to polynomial families such as $P_{k}(z)$,
$W_k(z)$ and $Q_k(z)$.  In this section we prove that $Q_{k}(z)$
and $W_{k}(z)$ have all their roots on the unit circle, by
extending the theorems used in the previous section. Let us
briefly recall how to derive results such as \eqref{Lakatos
condition 1} and \eqref{Lakatos condition 2}. For a real-valued
reciprocal polynomial $V_k(z)= \sum_{j=0}^k A_j z^j$, the
condition
\begin{align}\label{condition}
|A_{k}|>& \sum_{j=0}^{k}|cA_j -A_{k}|,
\end{align}
immediately implies that
\begin{equation}\label{Lakatos oscillating function}
1>\left|\frac{c}{A_k}V_k(z)-v_k(z)\right|,
\end{equation}
where $v_k(z)= \frac{z^{k+1}-1}{z-1}$.  Notice that if $v_{k}(z)$
is expanded in a geometric series, then \eqref{Lakatos oscillating
function} can be derived from \eqref{condition} as a simple
consequence of the triangle inequality. Despite the fact that
\eqref{Lakatos oscillating function} does not imply
\eqref{condition}, it turns out that \eqref{Lakatos oscillating
function} is \textit{easily} strong enough to conclude that
$V_{k}(z)$ has all of its zeros on the unit circle.  To
demonstrate this, first restrict $z$ to the unit circle, and write
$z=e^{i\theta}$ with $\theta \in (0, 2\pi)$, and
\[\tilde{v}_k(\theta)=z^{-(k+1)/2}v_k(z)= \frac{\sin\left(\frac{(k+1)\theta}{2} \right)}{\sin\left(\frac{\theta}{2} \right)}.\]
If $j<2k+2$ is a positive odd integer, then it is easy to show
that $\tilde{v}_k\left(\frac{j\pi}{k+1}\right)$ has sign
$(-1)^{(j-1)/2}$, and
$\left|\tilde{v}_k\left(\frac{j\pi}{k+1}\right)\right|\ge 1$. This
implies that $\tilde{v}_{k}(\theta)$ has at least $k+1$
interlacing positive and negative values in the interval
$(0,2\pi)$, and it has absolute value $\ge 1$ at each of those
points. By \eqref{Lakatos oscillating function} we can write
$\frac{c}{A_k}z^{-(k+1)/2}V_k(z)=\tilde{v}_k(\theta) +ET$, where
the error term has absolute value less than 1.  It follows that
$\frac{c}{A_k}z^{-(k+1)/2}V_k(z)$ changes sign at least $k$ times
for $\theta\in(0,2\pi)$. By the intermediate value theorem we
conclude that $V_{k}(z)$ has at least $k$ zeros on the unit
circle.  Since the polynomial has at most $k$ zeros, all of its
zeros must lie on the unit circle.

We can easily extend this idea by selecting a different
$v_k(z)$.\footnote{This principle was inspired by a careful study
of the proof in \cite{MSW}}  This typically entails constructing
$v_{k}(z)$ to approximate specific polynomial families.

\begin{definition}
Let $f(\theta) : (a,b)\rightarrow \mathbb{R}$ be a continuous
function. We call $f(\theta)$ a {\em $k$th order alternating
function on $(a,b)$}, if it assumes alternating positive and
negative (or negative and positive) values at points $p_i$, where
$a<p_1<\dots <p_{k+1}<b$. We say that $f(\theta)$ has
\textit{oscillation distance $d$}, if $|f(p_i)|>d$ for each
$i\in\{1,\dots,k+1\}$.
\end{definition}

\begin{lemma}\label{lemmaalternating} Suppose that $f(\theta)$ is a $k$th order alternating function on $(a,b)$, with oscillation distance $d$.
Let $F(\theta) : (a,b) \rightarrow \mathbb{R}$ be a continuous
function such that $|F(\theta)-f(\theta)|<d$ for all $\theta$.
Then $F(\theta)$ has at least $k$ zeros.
\end{lemma}
\begin{proof}  This lemma is essentially a restatement of the intermediate value theorem.  The proof follows immediately from the method described in the
previous discussion.
\end{proof}

\subsection{The zeros of $W_k(z)$ lie on the unit circle}

The main result of this subsection is the following theorem:
\begin{theorem}\label{Theorem on Wk}  Suppose that $k\ge 2$. The polynomial
\begin{align*}
W_{k}(z)=&(2^{2k-1}+2)P_k(z)-2^{2k}P_k(z/2)- P_k(2z)\\
=&\frac{(2\pi)^{2k-1}2^{2k}}{(2k)!}\sum_{j=0}^k (-1)^j
B_{2j}B_{2k-2j} (1-2^{1-2j})(1-2^{1-2k+2j})\binom{2k}{2j} z^{2j}
\end{align*}
has all of its zeros on the unit circle.  Furthermore, all of the
zeros are simple.
\end{theorem}

In order to prove Theorem \ref{Theorem on Wk} we first need to
establish that a certain trigonometric polynomial possesses the
alternating property with oscillation distance $0.3$.

\begin{lemma}\label{gk lemma} Suppose that $k>10$.  The function
\[w_k(\theta):= 2 \cos(k\theta) +\frac{\pi^2}{3}  \cos((k-2)\theta)+\frac{2}{(1-2^{1-2k})} \frac{\sin((k-3)\theta)}{\sin \theta}\]
is an alternating function of order $2k$ on $(-\pi,\pi)$, with
oscillation distance $0.3$.
\end{lemma}
\begin{proof}  We need to demonstrate that $|w_k(\theta)|>.3$ for $2k+1$
values of $\theta\in(-\pi,\pi)$.  We must also show that the sign
of $w_{k}(\theta)$ alternates over successive points in this set.
Since $w_k(\theta)$ is even with respect to $\theta$, and since
$w_{k}(0)>3$, we only need to demonstrate that there are an
additional $k$ such points in $(0,\pi)$. Suppose that $k>10$, let
$\alpha$ be defined by
\begin{equation}\label{alpha def}
\alpha=\frac{1}{\pi}\arccos\left(\frac{.3}{\frac{\pi^2}{3}-2}\right)=.42\dots,
\end{equation}
and let
\begin{equation}\label{j0 def}
j_0=[(k-1)\alpha]+1.
\end{equation}
We claim that $w_k(\theta)$ satisfies the necessary conditions on
the following set of points:
\begin{align*}
\emph{S}=&\left\{\frac{\pi}{k-1},\dots,\frac{(j_{0}-1)\pi}{k-1}\right\}
\cup\left\{\frac{(j_0-1/2)\pi}{k-1},\dots,\frac{(k-j_{0}-1/2)\pi}{k-1}\right\}\\
&\quad\cup\left\{\frac{(k-j_0)\pi}{k-1},\dots,\frac{(k-2)\pi}{k-1}\right\}\cup\left\{\frac{(k-(1-\epsilon))\pi}{k-1}\right\},
\end{align*}
where $\epsilon>0$ is sufficiently small. First expand
$w_{k}(\theta)$ using trigonometric identities
\begin{align*}
w_k(\theta)=&\cos((k-1)\theta)\cos(\theta)
\left(\frac{\pi^2}{3}+2-\frac{4}{1-2^{1-2k}}\right)\\
&+\sin((k-1)\theta)\left(\left(\frac{\pi^2}{3}-2-\frac{4}{1-2^{1-2k}}\right)
\sin \theta +\frac{2}{1-2^{1-2k}}\csc \theta\right).
\end{align*}
Notice that $\emph{S}$ (essentially) arises from cases where
either $\cos((k-1)\theta)=0$, or $\sin((k-1)\theta)=0$.

Begin by considering the cases where $\theta=\frac{j\pi}{k-1}$ and
$j\in\{1,\dots,k-2\}$. Then
\begin{equation}\label{gk integer case}
w_k\left(\frac{j\pi}{k-1}\right)= (-1)^{j}
\cos\left(\frac{j\pi}{k-1}\right)
\left(\frac{\pi^2}{3}+2-\frac{4}{1-2^{1-2k}}\right).
\end{equation}
In order to have $|w_k\left(\frac{j\pi}{k-1}\right)|>0.3$, we need
to restrict $j$ so that
\begin{align*}
\left|\cos\left(\frac{j\pi}{k-1}\right)\right|>
\frac{0.3}{\frac{\pi^2}{3}+2-\frac{4}{1-2^{1-2k}}}.
\end{align*}
In other words we must have
\[\frac{j}{k-1} \not \in (\alpha_k, 1-\alpha_k),\]
where
\begin{equation*}
\alpha_k=\frac{1}{\pi}\arccos\left(\frac{0.3}{\frac{\pi^2}{3}+2-\frac{4}{1-2^{1-2k}}}\right).
\end{equation*}
Since $k>10$ we have $(k-1)(\alpha_k-\alpha)\ll 1$. Therefore it
is sufficient that
\begin{equation*}
\frac{j}{k-1} \not \in (\alpha, 1-\alpha),
\end{equation*}
where $\alpha$ is defined in \eqref{alpha def}.  This implies that
$j\in\{1,2,\dots,j_{0}-1\}\cup\{k-j_{0},\dots,k-2\}$, with $j_{0}$
defined in \eqref{j0 def}.  If $j\in\{1,2,\dots,j_{0}-1\}$, then
by \eqref{gk integer case} $w_{k}\left(\frac{j\pi}{k-1}\right)$
has sign $(-1)^j$. If $j\in\{k-j_{0},\dots,k-2\}$ then the cosine
in \eqref{gk integer case} contributes an extra minus sign, and
$w_{k}\left(\frac{j\pi}{k-1}\right)$ has sign $(-1)^{j+1}$.

Now consider the case where $\theta=\frac{(j-1/2)\pi}{(k-1)}$ and
$j\in\{j_0,\dots,k-j_0\}$. By elementary properties of
trigonometric functions, $w_{k}(\theta)$ reduces to
\begin{equation}\label{gk at half integer rationals}
\begin{split}
w_k\left(\frac{(j-1/2)\pi}{k-1}\right)=(-1)^{j+1}\left(\left(\frac{\pi^2}{3}-2-\frac{4}{1-2^{1-2k}}\right)
\sin\left(\frac{(j-1/2)\pi}{k-1}\right)\right.\\
\left.+\frac{2}{1-2^{1-2k}}\csc\left(\frac{(j-1/2)\pi}{k-1}\right)\right).
\end{split}
\end{equation}
In order to place a lower bound on this expression, first choose
an interval $(\beta,1-\beta)$, which contains the set of rational
numbers $\{\frac{j_0-1/2}{k-1},\dots,\frac{k-j_0-1/2}{k-1}\}$.
This can be accomplished by selecting
\begin{equation*}
\begin{split}
\beta          &=\left\{\begin{array}{ll}
                      \alpha &\text{if $j_0>\alpha(k-1)+\frac{1}{2}$},\\
                      \alpha-\frac{1}{2(k-1)} &\text{if
                      $j_0<\alpha(k-1)+\frac{1}{2}$}.
                \end{array}\right.
\end{split}
\end{equation*}
Notice that one of these situations must occur, because \eqref{j0
def} guarantees that $j_0\in(\alpha(k-1),\alpha(k-1)+1)$. We
obtain the following lower bound from \eqref{gk at half integer
rationals}:
\begin{align*}
\left|w_k\left(\frac{(j-1/2)\pi}{k-1}\right)\right|\ge\min_{\theta\in(\pi\beta,\pi(1-\beta))}\left|\left(\frac{\pi^2}{3}-2-\frac{4}{1-2^{1-2k}}\right)
\sin\theta+\frac{2}{1-2^{1-2k}}\csc\theta\right|.
\end{align*}
The right-hand side is minimized at the end points of the interval
$(\pi \beta,\pi(1-\beta))$, so it follows that
\begin{align*}
\left|w_k\left(\frac{(j-1/2)\pi}{k-1}\right)\right|\ge\left|\left(\frac{\pi^2}{3}-2-\frac{4}{1-2^{1-2k}}\right)
\sin\pi\beta+\frac{2}{1-2^{1-2k}}\csc\pi\beta\right|
\end{align*}
Consider both choices of $\beta$, and recall the assumption that
$k>10$.  A few easy calculations are sufficient to obtain
\begin{equation*}
\begin{split}
\left|w_k\left(\frac{(j-1/2)\pi}{k-1}\right)\right|&>\left\{\begin{array}{ll}
                      .57 &\text{if $j_0>\alpha(k-1)+\frac{1}{2}$},\\
                      .34 &\text{if
                      $j_0<\alpha(k-1)+\frac{1}{2}$},
                \end{array}\right.
\end{split}
\end{equation*}
for all values of $j\in\{j_0,\dots,k-2j_0\}$.  It is easy to
deduce from \eqref{gk at half integer rationals} that the sign of
$w_{k}\left(\frac{(j-1/2)\pi}{k-1}\right)$ is precisely
$(-1)^{j}$.

Finally consider the value of
$w_k\left(\frac{(k-(1-\epsilon))\pi}{k-1}\right)$. Notice that
\begin{equation*}
w_{k}(\pi)=(-1)^k\left(\frac{\pi^2}{3}+2+\frac{2k-6}{1-2^{1-2k}}\right).
\end{equation*}
Since $k>10$ it follows easily that $|w_{k}(\pi)|>19$, and
$w_{k}(\pi)$ has sign $(-1)^k$.  If $\epsilon$ is sufficiently
small then $w_{k}\left(\frac{(k-(1-\epsilon))\pi}{k-1}\right)$
also has sign $(-1)^k$, and absolute value much larger than $.3$.

To briefly summarize the sign values of $w_k(\theta)$, we have the
following table:
\begin{equation*}
    \begin{tabular}{|c|c|p{6 in}|}
        \hline
         $\theta$ &   $Sign\left(w_k(\theta)\right)$\\
        \hline\hline
         $0$ &   $(-1)^{0}$\\
         $\frac{\pi}{k-1}$ &   $(-1)^{1}$\\
         $\frac{2\pi}{k-1}$ &   $(-1)^2$\\
         $\vdots$           & $\vdots$\\
         $\frac{(j_0-1)\pi}{k-1}$ &   $(-1)^{j_0-1}$\\
         $\frac{(j_0-1/2)\pi}{k-1}$ &   $(-1)^{j_0}$\\
         $\vdots$       &  $\vdots$\\
         $\frac{(k-j_0-1/2)\pi}{k-1}$ &   $(-1)^{k-j_0}$\\
         $\frac{(k-j_0)\pi}{k-1}$ &   $(-1)^{k-j_0+1}$\\
         $\vdots$& $\vdots$\\
         $\frac{(k-2)\pi}{k-1}$ &   $(-1)^{k-1}$\\
         $\frac{(k-(1-\epsilon))\pi}{k-1}$&$(-1)^k$\\
        \hline
    \end{tabular}
\end{equation*}
This table shows that $w_{k}(\theta)$ changes sign at least $k$
times over the interval $(0,\pi)$.
\end{proof}

Now we use Lemma \ref{gk lemma} to prove our main result.

\begin{proof}{(Theorem \ref{Theorem on Wk})}. Let us define $A_j$ using
\begin{align*}
W_{k}(i
z)&=\frac{(2\pi)^{2k-1}2^{2k}}{(2k)!}\sum_{j=0}^{k}B_{2j}B_{2k-2j}
(1-2^{1-2j})(1-2^{1-2k+2j}){2k\choose 2j}z^{2j}\\
&=\sum_{j=0}^{k}A_j z^{2j}.
\end{align*}
By  Lemma \ref{lemmaalternating} it suffices to prove that
$\left|\frac{z^{-k}W_k(iz)}{A_0}-w_k(z)\right|<0.3$ where
\begin{equation}\label{eq:h}
w_k(z)=(z^{k}+z^{-k})+ \frac{\pi^2}{6} (z^{k-2}+z^{2-k}) +
\frac{2}{(1-2^{1-2k})}  \frac{z^{k-3}-z^{3-k}}{z-z^{-1}},
\end{equation}
and $z=e^{i\theta}$.  Thus we write
\begin{align}
\left|\frac{z^{-k}W_k(iz)}{A_0}-w_k(z)\right|=&\left|\sum_{j=0}^k
\frac{A_j}{A_0} z^{2j-k}
 -(z^{k}+z^{-k})-\frac{\pi^2}{6} (z^{k-2}+z^{2-k})\right.\notag\\
 &\quad \left.-
\frac{2}{(1-2^{1-2k})} \frac{z^{k-3}-z^{3-k}}{z-z^{-1}}\right|\notag\\
&\leq 2 \left|
\frac{A_1}{A_0}-\frac{\pi^2}{6}\right|+\sum_{j=2}^{k-2}
\left|\frac{A_j}{A_0} - \frac{2}{(1-2^{1-2k})} \right| ,\label{Wgk
upper bound}
\end{align}
where the second step makes use of a geometric series and the
triangle inequality.

If we recall that $B_{0}=-1/2$, and use both inequalities from
\eqref{ineq:bernoulli}, then we find
\[\frac{A_j}{A_0} =\frac{B_{2j}B_{2k-2j}(1-2^{1-2j})(1-2^{1-2k+2j})\binom{2k}{2j}}{-B_{2k}(1-2^{1-2k})}< \frac{2}{(1-2^{1-2k})}. \]
Additionally we have, by \eqref{zeta},
\[\left| \frac{A_1}{A_0}-\frac{\pi^2}{6}\right|=\left|\frac{\pi^2}{6}\frac{\zeta(2k-2)}{\zeta(2k)}-\frac{\pi^2}{6} \right|.\]
This second expression goes to zero as $k\rightarrow \infty$. For
example, it is not hard to see that the absolute value is less
than $0.01$ for $k>4$.

Therefore we can remove the absolute value signs from \eqref{Wgk
upper bound}. We find that
\begin{align*}
\left|\frac{z^{-k}W_k(iz)}{A_0}-w_k(z)\right| \leq & 2
\frac{A_1}{A_0}-\frac{\pi^2}{3}+\frac{2(k-3)}{1-2^{1-2k}} -\sum_{j=2}^{k-2} \frac{A_j}{A_0}\\
=&2+4\frac{A_1}{A_0} -\frac{\pi^2}{3}+\frac{2(k-3)}{1-2^{1-2k}} -\sum_{j=0}^{k} \frac{A_j}{A_0}\\
=&2+
\frac{2\pi^2}{3}\frac{\zeta(2k-2)}{\zeta(2k)}-\frac{\pi^2}{3}+\frac{2(k-3)}{1-2^{1-2k}} -\frac{2k-1}{1-2^{1-2k}}\\
 =&2+\frac{2\pi^2}{3}\frac{\zeta(2k-2)}{\zeta(2k)}-\frac{\pi^2}{3}-\frac{5}{1-2^{1-2k}}\\
\leq& -3 +\left(2
\frac{1-2^{-2k}}{1-2^{3-2k}}-1\right)\frac{\pi^2}{3}.
\end{align*}
Notice that we evaluated $\sum_{j}A_j$ using the same Bernoulli
convolution identity \eqref{bernoulliconvolution}.  In addition, we have used the inequality
\[ \frac{1}{1-2^{-n}} < \zeta(n)<\frac{1}{1-2^{1-n}},\]
which are easy to deduce from the Euler product formula and the Dirichlet eta function. 
As
$k\rightarrow \infty$ this final upper bound approaches a limit of
$\frac{\pi^2}{3}-3\approx.289$. It is easy to see that it becomes
$<0.3$ for $k>6$.

In summary, we have proved the theorem for $k>10$. The cases for
$k\le 10$ are easily checked with the numerical method outlined in
section \ref{Section: partial Pk results}.
\end{proof}

\subsection{The zeros of $Q_k(z)$ lie on the unit circle}

Notice that the coefficients of $Q_{k}(z)$ involve the odd values
of the Riemann zeta function.  The primary goal of this subsection
is to prove the following theorem:

\begin{theorem}\label{Theorem on Qk}  Suppose that $k\ge 2$. The polynomial
\begin{align*}
Q_{k}(z)=&(2^{2k}+1)P_k(z)-2^{2k}P_k(z/2)- P_k(2z)\\
=&\frac{(2\pi)^{2k-1}}{(2k)!}\sum_{j=0}^k (-1)^j B_{2j}B_{2k-2j}
(2^{2j}-1)(2^{2k-2j}-1)\binom{2k}{2j} z^{2j}\\
&+ \zeta(2k-1)(2^{2k-1}-1)((-1)^k z+z^{2k-1}),
\end{align*}
has all of its non-zero roots on the unit circle.  Furthermore,
all of the zeros are simple.
\end{theorem}

As in the proof of Theorem \ref{Theorem on Wk}, the first step is
to construct a trigonometric polynomial which approximates
$Q_{k}(z)$.  Notice that $Q_{k}(z)$ has degree $2k-1$, and that it
has a trivial zero at $z=0$.  Therefore we need to prove that it
has $2k-2$ zeros on the unit circle.

\begin{lemma}\label{hk lemma} Suppose that $k>5$.  Then
\begin{equation*}
q_k(\theta):= 2 \cos((k-2)\theta) +\frac{4}{\pi}
\sin((k-1)\theta)+\frac{8(1-2^{3-2k})}{\pi^2(1-2^{2-2k})}\frac{\sin((k-3)\theta)}{\sin
\theta}
\end{equation*}
is an alternating function of order $2k-2$ on $(-\pi,\pi)$, with
oscillation distance $0.03$.
\end{lemma}
\begin{proof}
We need to demonstrate that $|q_k(\theta)|>.03$ for $2k-1$ values
of $\theta\in(-\pi,\pi)$.  We must also show that the sign of
$q_{k}(\theta)$ alternates over successive points in this subset.
The proof is similar to the proof of Lemma 2, so we will be brief.
Let $\alpha$ be defined by
\begin{equation}\label{alpha hk def}
\alpha=\frac{1}{\pi}\arccos\left(\frac{.03}{2-\frac{16}{\pi^2}}\right)=.47\dots,
\end{equation}
and let
\begin{equation}\label{j0 hk def}
j_0=[(k-1)\alpha]+1.
\end{equation}
We claim that $|q_k(\theta)|>.03$ on the following set of $2k+1$
points:
\begin{align*}
\emph{S}=&\left\{0\right\}\cup\left\{\pm\frac{\pi}{k-1},\dots,\pm\frac{(j_{0}-1)\pi}{k-1}\right\}
\cup\left\{\pm\frac{(j_0-1/2)\pi}{k-1},\dots,\pm\frac{(k-j_{0}-1/2)\pi}{k-1}\right\}\\
&\quad\cup\left\{\pm\frac{(k-j_0)\pi}{k-1},\dots,\pm\frac{(k-2)\pi}{k-1}\right\}\cup\left\{\pm\frac{(k-(1-\epsilon))\pi}{k-1}\right\}.
\end{align*}
If we consider
$S\setminus\{\frac{(j_0-1/2)\pi}{k-1},\frac{(k-j_0-1/2)\pi}{k-1}\}$,
then we obtain a subset of $2k-1$ points where the sign of
$q_{k}(\theta)$ alternates over successive points.

In order to prove this claim, first expand $q_k(\theta)$ using
trigonometric identities
\begin{equation*}
\begin{split}
 q_k(\theta)=&2 \cos((k-1)\theta)\cos(\theta) \left(1
-\frac{8(1-2^{3-2k})}{\pi^2(1-2^{2-2k})}\right)\\
&+\sin((k-1)\theta)\left(2\sin(\theta)\left(1-\frac{8(1-2^{3-2k})}{\pi^2(1-2^{2-2k})}
\right)  +\frac{4}{\pi}
+\frac{8(1-2^{3-2k})}{\pi^2(1-2^{2-2k})}\csc\theta \right).
\end{split}
\end{equation*}
Now consider the cases where $\theta=\frac{j\pi}{k-1}$, with
$-(k-2)\le j\le(k-2)$ and $j\ne0$. We have
\begin{equation}\label{hk integer points}
q_k\left(\frac{j\pi}{k-1}\right)= 2 (-1)^j
\cos\left(\frac{j\pi}{k-1}\right) \left(1
-\frac{8(1-2^{3-2k})}{\pi^2(1-2^{2-2k})}\right).
\end{equation}
To ensure that $|q_k\left(\frac{j\pi}{k-1}\right)|>0.03$, we need
to restrict $j$ so that
\[\left|\cos\left(\frac{j\pi}{k-1}\right)\right|
>\frac{0.03}{2\left(1
-\frac{8(1-2^{3-2k})}{\pi^2(1-2^{2-2k})}\right)}.\]
 By similar reasoning as in the proof of Lemma \ref{gk lemma}, it
is sufficient that
\[\frac{j}{k-1} \not \in (-(1-\alpha),-\alpha)\cup (\alpha, 1-\alpha) ,\]
where $\alpha$ is defined in \eqref{alpha hk def}.  This
immediately implies that $j\in\{\pm 1,\dots,\pm (j_0-1)\}\cup\{\pm
(k-j_0),\dots, \pm (k-2)\}$, where $j_0$ is defined in \eqref{j0
hk def}.  A careful inspection of \eqref{hk integer points}
reveals that the function has sign $(-1)^j$ for $j\in\{\pm
1,\dots,\pm(j_0-1)\}$, and sign $(-1)^{j+1}$ for
$j\in\{\pm(k-j_0),\dots,\pm (k-2)\}$.

Next consider the cases where $\theta=\pm\frac{(j-1/2)\pi}{(k-1)}$
and $j\in\{j_0,\dots,(k-j_0)\}$. We obtain
\begin{equation}\label{hk odd cases}
\begin{split}
q_k\left(\pm\frac{(j-1/2)\pi}{2(k-1)}\right)=(-1)^{j-1}&\left(\pm
2\sin\left(\frac{(j-1/2)\pi}{(k-1)}\right)\left(1-\frac{8(1-2^{3-2k})}{\pi^2(1-2^{2-2k})}
\right)\right.\\ &\quad\left.
\pm\frac{8(1-2^{3-2k})}{\pi^2(1-2^{2-2k})}\csc\left(\frac{(j-1/2)\pi}{(k-1)}
\right)+\frac{4}{\pi}\right).
\end{split}
\end{equation}
Now select $\beta$ so that
$\left\{\frac{j_0-1/2}{k-1},\dots,\frac{k-j_0-1/2}{k-1}\right\}\subset(\beta,1-\beta)$.
 Following Lemma \ref{gk lemma}, this is accomplished by selecting
\begin{equation*}
\begin{split}
\beta          &=\left\{\begin{array}{ll}
                      \alpha &\text{if $j_0>\alpha(k-1)+\frac{1}{2}$},\\
                      \alpha-\frac{1}{2(k-1)} &\text{if
                      $j_0<\alpha(k-1)+\frac{1}{2}$}.
                \end{array}\right.
\end{split}
\end{equation*}
Therefore we obtain
\begin{equation*}
\begin{split}
\left|q_k\left(\pm\frac{(j-1/2)\pi}{2(k-1)}\right)\right|&\ge\min_{\theta\in(\pi\beta,\pi(1-\beta))}\left|2\left(1-\frac{8(1-2^{3-2k})}{\pi^2(1-2^{2-2k})}
\right)\sin\theta
+\frac{8(1-2^{3-2k})}{\pi^2(1-2^{2-2k})}\csc\theta \pm\frac{4}{\pi}\right|\\
&\ge\left|2\left(1-\frac{8(1-2^{3-2k})}{\pi^2(1-2^{2-2k})}
\right)\sin\pi\beta
+\frac{8(1-2^{3-2k})}{\pi^2(1-2^{2-2k})}\csc\pi\beta \pm
\frac{4}{\pi}\right|.
\end{split}
\end{equation*}
Checking both possible values of $\beta$, and both possible signs
of $\pm$, leads to a lower bound which holds for $k>5$:
\[\left|q_k\left(\pm\frac{(2j-1)\pi}{2(k-1)}\right)\right|>
0.08...\] The final signs are summarized in the table below.

The only remaining cases are when $j\in\{0\}\cup\{\pm
\frac{(k-(1-\epsilon))\pi}{k-1}\}$. These cases can be easily
dispensed with by elementary properties of trigonometric
functions.

To briefly summarize the sign values of $q_k(\theta)$, we have the
following table:
\begin{equation*}
    \begin{tabular}{|c|c|c|c|p{6 in}|}
        \hline
         $\theta$ &   $Sign\left(q_k(\theta)\right)$&$\theta$ &   $Sign\left(q_k(\theta)\right)$\\
        \hline\hline
        $-\frac{(k-(1-\epsilon))\pi}{k-1}$ &$(-1)^k$ &$0$ &   $(-1)^{0}$\\
         $-\frac{(k-2)\pi}{k-1}$& $(-1)^{k-1}$ &$\frac{\pi}{k-1}$ &   $(-1)^{1}$\\
        $\vdots$&$\vdots$ &$\frac{2\pi}{k-1}$ &   $(-1)^2$\\
         $-\frac{(k-j_0)\pi}{k-1}$& $(-1)^{k-j_0+1}$ & $\vdots$           & $\vdots$\\
         $-\frac{(k-j_0-1/2)\pi}{k-1}$& $(-1)^{k-j_0}$ &$\frac{(j_0-1)\pi}{k-1}$ &   $(-1)^{j_0-1}$\\
         $\vdots$ & $\vdots$ &$\frac{(j_0-1/2)\pi}{k-1}$ &   $(-1)^{j_0-1}$\\
        $ -\frac{(j_0-1/2)\pi}{k-1}$& $(-1)^{j_0}$&$\vdots$       &  $\vdots$\\
         $-\frac{(j_0-1)\pi}{k-1}$& $(-1)^{j_0-1}$ & $\frac{(k-j_0-1/2)\pi}{k-1}$ &   $(-1)^{k-j_0-1}$\\
         $\vdots$&$\vdots$  &$\frac{(k-j_0)\pi}{k-1}$ &   $(-1)^{k-j_0+1}$\\
         $-\frac{2\pi}{k-1}$&$(-1)^{2}$ &$\vdots$& $\vdots$\\
         $-\frac{\pi}{k-1}$&$(-1)^{1}$ &$\frac{(k-2)\pi}{k-1}$ &   $(-1)^{k-1}$\\
         & &$\frac{(k-(1-\epsilon))\pi}{k-1}$&$(-1)^k$\\
        \hline
    \end{tabular}
\end{equation*}
Notice that there are precisely $2k+1$ values of $\theta$ in this
table.  If we exclude the cases where
$\theta\in\{\frac{(j_0-1/2)\pi}{k-1},\frac{(k-j_0-1/2)\pi}{k-1}\}$,
then the sine of $q_k(\theta)$ alternates over the remaining
$2k-1$ values of $\theta$.
\end{proof}

Next we use Lemma \ref{hk lemma} to establish that $Q_k(z)$ has
all of its non-zero roots on the unit circle for $k\ge2$.

\begin{proof}(Theorem \ref{Theorem on Qk}).
Let us define $A_j$ using
\begin{align*}
Q_{k}(i z) =&\frac{(2\pi)^{2k-1}}{(2k)!}\sum_{j=0}^k
B_{2j}B_{2k-2j}
(2^{2j}-1)(2^{2k-2j}-1)\binom{2k}{2j} z^{2j}\\
&+i(-1)^k \zeta(2k-1)(2^{2k-1}-1)(z-z^{2k-1})\\
=&\sum_{j=0}^{k}A_j z^{2j}+i(-1)^k
\zeta(2k-1)(2^{2k-1}-1)(z-z^{2k-1}).
\end{align*}
In order to simplify the following analysis, we have intentionally
defined $A_j$ to only involve the even coefficients of $Q_k(i z)$.
Notice that $A_0=A_{k}=0$, and that
\begin{align*}
A_1=& \frac{(2\pi)^{2k-1}}{(2k)!} B_{2k-2} (2^{2k-2}-1)
\frac{k(2k-1)}{2}\\
 =& (-1)^k\zeta(2k-2)(2^{2k-2}-1).
\end{align*}
 Suppose that $k>5$.  Then by Lemma
\ref{lemmaalternating} it suffices to prove that
$\left|\frac{z^{-k}Q_k(iz)}{A_1}-q_k(z)\right|<0.03$ where
\begin{equation}\label{eq:h}
q_k(z)=(z^{k-2}+z^{2-k})- \frac{2i}{\pi} (z^{k-1}-z^{1-k})
+\frac{8}{\pi^2}
\frac{(1-2^{3-2k})}{(1-2^{2-2k})}\frac{z^{k-3}-z^{3-k}}{z-z^{-1}},
\end{equation}
and $z=e^{i\theta}$.

Therefore let us begin by writing
\begin{align*}
\left|\frac{z^{-k}Q_z(iz)}{A_1}-q_k(z)\right|=&
\left|\sum_{j=0}^k \frac{A_j}{A_1} z^{2j-k}+(-1)^k i \frac{\zeta(2k-1)(2^{2k-1}-1)}{A_1}(z^{1-k}- z^{k-1})\right.\\
& \left.- (z^{k-2}+z^{2-k})+ \frac{2i}{\pi} (z^{k-1}-z^{1-k})-
\frac{8}{\pi^2}\frac{(1-2^{3-2k})}{(1-2^{2-2k})}\frac{z^{k-3}-z^{3-k}}{z-z^{-1}}\right|\\
&\leq \sum_{j=2}^{k-2} \left|\frac{A_j}{A_1} -
\frac{8}{\pi^2}\frac{(1-2^{3-2k})}{(1-2^{2-2k})}\right| + 2 \left|
(-1)^k\frac{\zeta(2k-1)(2^{2k-1}-1)}{A_1}-\frac{2}{\pi}\right|.
\end{align*}
The second step follows from substituting a geometric series, and
then applying the triangle inequality. We know from equation
\eqref{Thm2 Aj Ak inequality} (after noting the change in
definition of $A_j$), that
\[\frac{A_j}{A_1} > \frac{8}{\pi^2}\frac{(1-2^{3-2k})}{(1-2^{2-2k})}.\]
In addition
\[\left| (-1)^k\frac{\zeta(2k-1)(2^{2k-1}-1)}{A_1}-\frac{2}{\pi}\right|
= \left|
\frac{2}{\pi}\frac{\zeta(2k-1)(1-2^{1-2k})}{\zeta(2k-2)(1-2^{2-2k})}-\frac{2}{\pi}\right|.\]
This limit tends to zero. A simple calculation shows that the
quantity is less than $0.001$ for $k>3$.

Therefore we can remove the absolute value signs from the
inequality.  We are left with
\begin{align*}
&\left|\frac{z^{-k}Q_z(iz)}{A_1}-q_k(z)\right|\\
 &\qquad\leq
\sum_{j=2}^{k-2} \frac{A_j}{A_1} -
(k-3)\frac{8}{\pi^2}\frac{(1-2^{3-2k})}{(1-2^{2-2k})}
+\frac{4}{\pi} -
2(-1)^k\frac{\zeta(2k-1)(2^{2k-1}-1)}{A_1}\\
&\qquad=\sum_{j=1}^{k-1} \frac{A_j}{A_1} -2 -
(k-3)\frac{8}{\pi^2}\frac{(1-2^{3-2k})}{(1-2^{2-2k})}
+\frac{4}{\pi}-\frac{4}{\pi}\frac{\zeta(2k-1)(1-2^{1-2k})
}{\zeta(2k-2)(1-2^{2-2k})}\\
&\qquad=(2k-1)\frac{4}{\pi^2}\frac{\zeta(2k)(1-2^{-2k})}{\zeta(2k-2)(1-2^{2-2k})}-2
- (k-3)\frac{8(1-2^{3-2k})}{\pi^2(1-2^{2-2k})}\\
&\qquad\quad+\frac{4}{\pi}-\frac{4}{\pi}\frac{\zeta(2k-1)(1-2^{1-2k})}{\zeta(2k-2)(1-2^{2-2k})}.
\end{align*}
As usual, we have evaluated $\sum_{j}A_j$ using
\eqref{bernoulliconvolution}.  The limit of the upper bound is
$\frac{20}{\pi^2}-2$. It is easy to see that it becomes $<0.03$
for $k>6$. The cases for $k<6$ are easily proved with the
numerical method described in Section \ref{Section: partial Pk
results}.
\end{proof}

We conclude this section by deriving a second approximation for
$\zeta(3)/\pi^3$.  If we begin with the expression for $Q$ in
terms of $P$, and then substitute a truncated version of
\eqref{zeta(3) series} for $P$, we can obtain
\begin{equation}\label{zeta 3 approx 2}
\zeta(3)\approx\frac{z}{1+z^2}\frac{\pi^3}{14},
\end{equation}
where
\begin{equation}
0=\frac{2 z}{e^{4 \pi z}-1} + \frac{8 z^3}{e^{\pi/z}-1} -
 \frac{17z}{e^{2 \pi z}-1} -\frac{17 z^3}{e^{2 \pi/z}-1} +
 \frac{8 z}{ e^{\pi z}-1} + \frac{2 z^3}{e^{4 \pi/z}-1}.
\end{equation}
Selecting the zero given by $z\approx .92-.39 i$, yields $4$
decimal places of accuracy in \eqref{zeta 3 approx 2}. Notice that
this approximation is slightly worse than \eqref{zeta 3 approx 1}.

\section{Partial results on $P_{k}(z)$}\label{Section: partial Pk results}

We have made a number of unsuccessful attempts to apply the
theorems of Schinzel, Lakatos and Losonczi, and their
generalizations to the case of $P_{k}(z)$.\footnote{The most
general result that we can prove is that $P_k(z)$ has at least
$k-1$ zeros in half of the unit circle, using the construction from \cite{MSW}.} A piece of evidence
indicating that these methods may not work is given by the result
in \cite{MSW} which shows that $R_{2k+1}(z)$ has four zeros that do
not lie in the unit circle (by comparison $Y_k(z)$ has all of its
roots on the unit circle).

We will briefly describe one instance were Lakatos's condition
\eqref{Lakatos condition 1} fails, because it leads to an
interesting formula. Let us define the $A_j$'s as follows:
\begin{equation*}
\begin{split}
\sum_{j=0}^{4k}A_j z^j&=\left|P_{k}(i z)\right|^2\\
&=\left(\frac{2}{\pi}\sum_{j=0}^{k}\zeta(2j)\zeta(2k-2j)
z^{2j}\right)^2+\zeta^2(2k-1)\left(z^{2k-1}-z\right)^2.
\end{split}
\end{equation*}
Notice that $|P_k(i z)|^2$ has all of its zeros on the unit circle,
if and only if $P_k(z)$ also has all of its zeros on the unit
circle.  Computational experiments helped us to make the following
observation:
\begin{observation}Suppose that $k\ge 2$, then
\begin{equation}\label{Aj conjecture for P}
4k(k-1)|A_{4k}|=\sum_{j=0}^{4k}\left|A_{4k}-A_{j}\right|.
\end{equation}
\end{observation}
Formula \eqref{Aj conjecture for P} can be proved with the
identities for Bernoulli numbers that we used in Theorem
\ref{Theorem on qk}. This immediately rules out the possibility of
applying \eqref{Lakatos condition 1} (condition \eqref{Lakatos
condition 2} can also be ruled out by slightly different methods).
It is curious to note that the right-hand side of \eqref{Aj
conjecture for P} appears to involve odd zeta values, whereas the
left-hand side does not.  It turns out that when \eqref{Aj
conjecture for P} is explicitly calculated, the odd zeta values
drop out.

\begin{theorem}\label{Pk partial theorem} Suppose that $2\le k<1000$.  Then all of the zeros of $P_{k}(z)$
lie on the unit circle.  Furthermore, all of the zeros are simple.
\end{theorem}

While we have not been able to prove a general theorem concerning
$P_{k}(z)$, we have been able to prove Theorem \ref{Pk partial
theorem} for $k<1000$. The proof uses a standard computational
method based on the intermediate value theorem. Notice that the
map
\begin{equation*}
z+z^{-1}\rightarrow 2u,
\end{equation*}
sends the unit circle to the real interval $[-1,1]$.  Under this
transformation, we also have
\begin{equation*}
z^{k}+z^{-k}\rightarrow 2T_{k}\left(u\right),
\end{equation*}
where $T_k(u)$ is the usual Chebyshev Polynomial. If we write
$(z^{2k}+(-1)^k
)P_{k}(z)=z^{2k}\left(P_{k}(z)+P_{k}(1/z)\right)=2z^{2k}P_{k}^{*}(u)$,
then it follows that $P_{k}(z)$ has all of its zeros on the unit
circle, if and only if
\begin{equation*}
\begin{split}
P_k^{*}(u):=&\frac{(2\pi)^{2k-1}}{(2k)!}\sum_{j=0}^{k}(-1)^j
B_{2j}B_{2k-2j}{2k\choose
2j}T_{2j}(u)\\
&+\zeta(2k-1)\left(T_{2k-1}(u)+(-1)^k T_{1}(u)\right),
\end{split}
\end{equation*}
has all of its zeros in the interval $[-1,1]$.  It is easy to
count real zeros of real-valued polynomials. The intermediate
value theorem allows one to find zeros by detecting sign changes.
Since $P_{k}^{*}(u)$ has degree $2k$, it is only necessary to
detect $2k$ sign changes in $[-1,1]$ (fewer sign changes are
required if zeros lie at $u=\pm 1$). We have successfully carried
out these calculations for $k<1000$.


\section{Conclusion}\label{section: Conclusion}


In conclusion, we have shown that $S_{k}(z)$, $Y_{k}(z)$,
$W_{k}(z)$, and $Q_{k}(z)$ have all of their non-zero roots on
the unit circle. These polynomials have a strong connection to the
Ramanujan polynomials. We were disappointed that we were unable to
to deduce a similar theorem concerning $P_{k}(z)$, however we are
hopeful that the approach outlined in Section \ref{Section on Wk
and Qk} might eventually succeed in this case. An additional
avenue might involve studying the zeros of a truncated version of
the right-hand side of \eqref{zeta(3) series}. Notice that the
zeros of $P_{k}(z)$ are very well approximated by the zeros of
\begin{equation*}
0=\frac{z^{k-1}}{e^{2\pi /z}-1}+(-1)^{k+1}\frac{ z^{1-k}}{e^{2\pi
z}-1}.
\end{equation*}
Thus, it should be a worthwhile endeavor to study the zeros of
these auxiliary functions.

\textit{Acknowledgements} The authors thank Wadim Zudilin for the
interesting discussions, encouragement and helpful comments, and M. Ram Murty for calling our attention to \cite{GMR}.  A portion of this
research was carried out while Mat Rogers was visiting the
University of Georgia.  He is grateful for their hospitality.

\end{document}